\documentclass[reqno,12pt]{amsart}

\usepackage{amsmath,amsthm,amssymb,amsfonts,epsfig,color,graphicx,enumerate,stackrel}
\usepackage[a4paper,margin=3.25truecm]{geometry}
\usepackage[T1]{fontenc}

\DeclareMathOperator{\dist}{dist}

\def\R{\mathbb{R}}

\def\O{\Omega}

\def\B{\mathcal{B}}

\def\HH{\mathcal{H}}

\def\LL{\mathcal{L}}

\def\V{\mathcal{V}}

\newcommand{\be}{\begin{equation}}
\newcommand{\ee}{\end{equation}}
\newcommand{\bib}[4]{\bibitem{#1}{\sc#2: }{\it#3. }{#4.}}

\numberwithin{equation}{section}
\theoremstyle{plain}
\newtheorem{theo}{Theorem}[section]

\theoremstyle{remark}
\newtheorem{rema}[theo]{\bf Remark}

\title[The problem of minimal resistance, old and new]{The problem of minimal resistance, old and new}

\author{Giuseppe Buttazzo}

\begin{document}

\maketitle

{\it Dedicated to Alexander Plakhov for his 65th birthday}

\begin{abstract}

Since its original formulation by Isaac Newton in 1685, the problem of determining bodies of minimal resistance moving through a fluid has been one of the classical problems in the calculus of variations. Initially posed for cylindrically symmetric bodies, the problem was later extended to general convex shapes, as explored in \cite{BK93}, \cite{BFK95}. Since then, this broader formulation has inspired a number of articles dedicated to the study of the geometric and analytical properties of optimal shapes, with particular attention to their structure, regularity, and behavior under various constraints. In this article, we provide a comprehensive overview of the principal results that have been established, highlighting the main theoretical advancements. Furthermore, we introduce some new directions of research, some of which were described in \cite{P12}, that offer promising perspectives for future investigation.

\end{abstract}



%%%%%%%%%%%%%%%%%%%%%%%%%%%%%%%%%%%%%%%%%%%%%%%%%%
\section{Introduction}\label{sintro}

The problem to determine the optimal shape of a body that minimizes resistance, whether aerodynamic in air or hydrodynamic in water, represents one of the earliest and most classical problems in the calculus of variations. This classical inquiry traces back to the seminal work of Sir Isaac Newton, who in 1685 provided a pioneering treatment of the problem. In his {\it Philosophiae Naturalis Principia Mathematica}, Newton introduced a remarkably simple yet insightful model for computing the resistance experienced by a body moving through an idealized medium, specifically, an inviscid (non-viscous) and incompressible fluid.

In his own words:\\
\smallskip
{\it ``If in a rare medium, consisting of equal particles freely disposed at equal distances from each other, a globe and a cylinder described on equal diameter move with equal velocities in the direction of the axis of the cylinder, (then) the resistance of the globe will be half as great as that of the cylinder. \dots I reckon that this proposition will be not without application in the building of ships.''}
\smallskip

This early formulation of the minimal resistance problem laid the groundwork for centuries of mathematical investigation and practical engineering application. The historical evolution of the problem is thoroughly described in works such as the book by Goldstine \cite{gold}, and its essence may be described in the following terms:

Given a body moving at constant velocity through a fluid medium, and assuming that the body's rear end is constrained to have a prescribed maximal cross-sectional area perpendicular to the direction of motion, determine the shape of the body that experiences the least possible resistance.

Naturally, the specific solution to this problem is deeply influenced by the manner in which resistance is defined. In order to render the problem mathematically tractable, Newton adopted a number of idealizing assumptions, including:

\begin{itemize}
\item[-] The fluid is modeled as consisting of discrete particles that are {\it mutually non-interacting}, that is each particle moves independently of the others. These particles travel in straight lines at a uniform velocity, all aligned parallel to the direction of the undisturbed stream. Equivalently, one may adopt a reference frame in which the fluid particles are considered to be at rest, uniformly distributed throughout space, while the solid body moves through this stationary medium with constant velocity. In this perspective, the motion of the body relative to the fluid generates the same sequence of interactions as in the original formulation, but the problem is analyzed from the vantage point of the moving object traversing a quiescent field of particles.

\item[-] The resistance encountered by the body arises solely from direct collisions, namely shocks between the incoming fluid particles and the surface of the body. These interactions are governed by the classical laws of {\it perfectly elastic} collisions, meaning that the kinetic energy and momentum of each particle are conserved in the interaction, with no energy lost to heat or deformation.

\item[-] All secondary and more complex physical phenomena, such as {\it tangential friction} (viscous shear stress), the generation of {\it vorticity}, and the onset of {\it turbulence}, are entirely neglected. The model thus omits any rotational or viscous effects in the flow.
\end{itemize}

While these assumptions render Newton's model a significant simplification, and indeed a rather crude approximation of real fluid dynamics, it nevertheless proves remarkably insightful in certain regimes. Specifically, it yields useful predictions in scenarios such as:
\begin{itemize}
\item[-] the motion of a body through a rarefied gas at relatively low speeds;
\item[-] the high-speed (supersonic or hypersonic) motion of bodies through an ideal, inviscid gas;
\item[-] the behavior of {\it slender bodies}, where the flow can be effectively approximated as one-dimensional along the principal axis of motion.
\end{itemize}

A comprehensive discussion of the physical validity and mathematical implications of these idealizations can be found in \cite{K02} and in the monograph by Miele \cite{miel}, which also presents several extensions and generalizations of the original Newtonian formulation. For further insight, particularly regarding applications to hypersonic aerodynamics, the reader is referred to the influential work of Hayes and Probstein \cite{haye}.

Under the simplifying assumptions described above, one can, through elementary trigonometric reasoning, derive what is known as the {\it Newtonian sine-squared pressure law}. This law occupies a central place in Newton's theory of minimal resistance, providing a quantitative expression for the pressure distribution on the surface of a body as it interacts with the surrounding medium.

Specifically, the law states that the pressure exerted by the fluid particles at a given point on the body's surface is directly proportional to the square of the sine of the angle of incidence. That is, the pressure at each point varies as $\sin^2\vartheta(x)$, where $\vartheta(x)$ denotes the local angle between the tangent to the surface (or, equivalently, the normal to the surface) and the direction $x$ of the undisturbed flow.

This simple and powerful relation encapsulates the geometric dependence of aerodynamic or hydrodynamic resistance in Newton's idealized setting, and forms the analytical foundation for identifying body profiles that minimize resistance in a non-interacting particle medium.

Let us now consider a more precise mathematical formulation of the problem under Newton's assumptions. Denote by $\Omega$ the prescribed maximal cross section of the body, assumed to lie in a horizontal plane and to be orthogonal to the direction of motion, which we take to be vertically upward. The front profile of the body, that is the portion first encountering the oncoming flow, can be described by a function $u:\O\to\R$, where $u(x)$ represents the vertical height of the body at the point $x\in \O$.

At any given point on the surface, specifically at $\big(x,u(x)\big)$, the incident particles collide with the surface and impart momentum to the body. Due to Newton's sine-squared pressure law, and assuming the impact of each particle occurs only once (that is no multiple collisions or reflections), one finds that the decelerating force exerted by such a collision is inversely proportional to the expression $1+|\nabla u(x)|^2)$.

Integrating this quantity over the entire cross-sectional domain $\O$, and accounting for the constant density $\rho$ of the fluid and the uniform velocity $v$ of the body, the total resistance force experienced by the body can be expressed as
$$\rho v^2\int_\O\frac{1}{1+|\nabla u(x)|^2}\,dx.$$
The problem of finding the optimal body shape that minimizes resistance then reduces to the variational problem
\begin{equation}\label{pbmin}
\min\bigg\{F(u)\ :\ u\text{ admissible}\bigg\},
\end{equation}
where $F(u)$ is the resistance functional
$$F(u)=\int_\O\frac{1}{1+|\nabla u(x)|^2}\,dx.$$

However, several nontrivial mathematical challenges arise in the analysis of this problem. Most notably, the integrand of the functional $F$ is neither convex nor coercive. As a result, the classical direct methods of the calculus of variations, typically relying on weak lower semicontinuity and coercivity to guarantee the existence of minimizers, are not directly applicable. This significantly complicates the search for an existence theorem.

The identification of an appropriate class of admissible functions for the minimization problem \eqref{pbmin} is a subtle and nontrivial matter. Without carefully imposed constraints, the problem becomes ill-posed and loses its physical and mathematical meaning.

To illustrate this, consider the scenario in which no restrictions are placed on the admissible functions $u$. In such a case, the minimization problem becomes degenerate. Indeed, for each positive integer $n$, let us define the sequence of functions
$$u_n(x)=n\dist(x,\partial\O),$$
where $\dist(x,\partial\O)$ denotes the Euclidean distance from the point $x\in\O$ to the boundary $\partial\O$. These functions represent increasingly steep conical profiles that peak in the interior of $\O$ and vanish along the boundary. A straightforward computation shows that
$$\lim_{n\to\infty}\int_\O\frac{1}{1+|\nabla u_n(x)|^2}\,dx=0,$$
while the resistance functional remains strictly positive for any nontrivial profile. This implies that the infimum of the resistance functional is zero, but this value cannot be attained by any admissible function $u$. Hence, no minimizer exists in this unrestricted setting, rendering the variational problem devoid of meaningful solutions.

To remedy this issue, it is natural to introduce a constraint on the height of the body's profile, such as bounding $u$ from above. A typical choice is to restrict the class of admissible functions to those satisfying
$$0\le u(x)\le M\qquad\text{for all }x\in\O,$$
for some fixed positive constant $M$. This condition enforces a physically reasonable bound on the maximal height of the body and ensures the geometry remains contained within a compact region.

However, even this constraint does not, by itself, guarantee the existence of a minimizer. Indeed, further complications arise when one considers sequences of oscillatory functions that remain bounded in height but oscillate with increasingly high frequency. As an illustrative example, consider the sequence
$$u_n(x)=M\sin^2(n|x|),$$
which clearly satisfies the height constraint $0\le u_n(x)\le M$. Yet, despite this uniform bound, we still have for the associated resistance functional
$$\lim_{n\to\infty}\int_\O\frac{1}{1+|\nabla u_n(x)|^2}\,dx=0,$$
indicating that the total resistance can be made arbitrarily small, even though the functions $u_n$ remain bounded between $0$ and $M$. In other words, no optimal profile exists in this setting, despite the boundedness constraint on the height.

In addition, the strong oscillations of the functions $u_n$ are in clear contradiction to Newton's fundamental physical assumption that each particle interacts with the body's surface at most once. Indeed, highly oscillating profiles give rise to multiple potential points of impact along the same trajectory, thereby violating the model's foundational premise. Thus, to obtain a well-posed variational problem with guaranteed existence of minimizers, it becomes necessary to impose additional structural conditions on the admissible profiles, so as to preclude such pathological oscillatory behavior and ensure some kind of compactness of minimizing sequences.

In the following sections we present the existence results of optimal shapes for the minimum resistance problem. However, Newton's model is flexible enough to allow for several generalizations, which are necessary when considering situations closer to real life. In particular, the following generalizations are of interest.

\begin{itemize}
\item[-] The optimization problem presented above is expressed within the framework of Cartesian coordinates, under the assumption that the boundary of the solid body $E$ can be represented as the graph of a suitably function $u$. While this formulation is convenient for certain analytical and computational purposes, it lacks invariance under geometric transformations and does not fully capture the intrinsic geometry of the problem. For a more natural and coordinate-independent perspective, it is desirable to reformulate the problem intrinsically, as a variational principle involving the minimization of a cost functional defined by an integral over the boundary $\partial E$ of the domain. This intrinsic formulation has been developed and rigorously analyzed in \cite{BG97}, and we shall briefly revisit it in what follows.

\item[-] In certain situations, the imposition of a prescribed height constraint is too restrictive, limiting the set of admissible configurations in a way that may not be physically realistic. As a result, alternative forms of constraints have been considered in the literature and are often more appropriate, depending on the specific nature of the problem under study. Among these, volume constraints, ensuring that the total measure of the admissible bodies $E$ remains fixed, or perimeter constraints, controlling the total surface area $|\partial E|$, are of particular interest. These types of constraints often arise naturally in physical models and variational formulations, and they allow for a richer and more flexible class of admissible competitors in the optimization problems.

\item[-] Beyond the basic {\it drag} minimization considered in the Newton's problem, there are other shape functionals of considerable importance in practical applications, particularly in the context of aerodynamics and fluid dynamics. One notable example is the {\it lift} generated by an airfoil or a solid body immersed in a flow, which plays a critical role in aircraft performance. In many design problems, while the primary objective may be to minimize the drag, it is essential to ensure that the lift remains above a certain prescribed threshold in order to maintain sufficient aerodynamic support.

Within the framework of Cartesian coordinates, and following the Newton geometrical approach, the lift generated along the $x_1$ axis, can be represented by the functional
$$L(u)=\int_\O\frac{\nabla u\cdot e_1}{1+|\nabla u|^2}\,dx,$$
where $e_1$ denotes the unit vector in the $x_1$ direction, and $u$ defines the shape via its graph. Under this formulation, the associated optimization problem takes the form
$$\min\bigg\{\int_\O\frac{1}{1+|\nabla u(x)|^2}\,dx\ :\ u\text{ admissible, }L(u)\ge L_0\bigg\},$$
where $L_0>0$ is a prescribed minimum lift requirement. This constrained minimization problem reflects the practical necessity of balancing aerodynamic efficiency (through drag reduction) with functional performance (via sufficient lift).

\item[-] In Newton's classical model, the surrounding fluid is idealized as being composed of stationary particles, with the motion attributed entirely to the solid body traversing through it. This idealization can be interpreted as a {\it zero-temperature} approximation, wherein thermal agitation of the fluid particles is entirely neglected. However, in realistic physical scenarios, the ambient fluid is at a {\it positive temperature}, and its constituent particles exhibit incessant, random motion. This thermal motion introduces additional effects that significantly influence the interaction between the fluid and the moving body. Consequently, Newton's model must be extended to accommodate this richer framework by incorporating the impact of temperature on the resistance functional. Such a generalization, accounting for thermal contributions to the dynamical behavior, has been rigorously developed in \cite{P04} and \cite{PT05}. A comprehensive exposition of this extended theory is also presented in the monograph \cite{P12}.

\item[-] Another significant limitation of Newton's model consists in its restrictive treatment of particle-body interactions, which produces the expression of the resistance functional. In the classical formulation, the resistance is computed solely based on the first collision between the fluid particles and the solid body, with no consideration given to subsequent interactions. While this assumption is entirely appropriate in the context of convex bodies, where multiple reflections are geometrically precluded, it becomes insufficient when the analysis is extended to broader classes of admissible shapes, including non-convex or even multiply connected bodies. In such cases, the possibility of multiple reflections or intricate scattering phenomena must be carefully accounted for, necessitating a profound revision of the original optimization framework. Moreover, this more general framework is connected to the study of {\it low visibility} phenomena in geometric optics, where one seeks configurations that minimize the detectability of a body by a stream of incoming particles. For a detailed presentation of these issues and their mathematical treatment, we refer the reader to the article \cite{P03} and to the comprehensive monograph \cite{P12}.

\end{itemize}

%%%%%%%%%%%%%%%%%%%%%%%%%%%%%%%%%%%%%%%%%%%%%%%%%%
\section{The Cartesian model}\label{scarte}

In this section, we turn our attention to the classical Newton's problem of minimal resistance, specifically for solid bodies whose exposed (or free) surface can be represented in Cartesian coordinates by
$$\Big\{y=u(x)\ :\ x\in\O\Big\}.$$
Here $\O\subset\R^d$ denotes a prescribed bounded, convex domain that describes the maximal cross-sectional area of the body. The function $u:\O\to\R$ describes the height profile of the body's surface above the domain $\O$, and it is assumed to be concave. This assumption of concavity is essential in the physical modeling of the problem, as it reflects the requirement that each fluid particle has at most a single impact with the body, a foundational principle in Newton's formulation.

Moreover, in order to ensure that the optimization problem is mathematically well-posed, it is necessary to impose a constraint on the maximal height that the body may attain. This upper bound is justified both by the considerations discussed in the introduction.

The optimal resistance problem is then
\be\label{mpcart}
\min\bigg\{\int_\O\frac{1}{1+|\nabla u|^2}\,dx\ :\ u\in C_M\bigg\},
\ee
where the class of admissible functions $C_M$ is defined by
\be\label{CM}
C_M=\Big\{u:\O\to\R,\text{ $u$ concave, }0\le u\le M\Big\}.
\ee
More generally, one may consider a broader class of variational problems in the same admissible set $C_M$, of the form
\be\label{mpphi}
\min\bigg\{\int_\O\phi(x,u,\nabla u)\,dx\ :\ u\in C_M\bigg\},
\ee
where $\phi:\R^d\to\R$ is a given integrand, possibly modeling more general physical effects or cost criteria beyond Newton's classical formulation.

A crucial feature of the admissible class $C_M$ is its compactness in suitable Sobolev function spaces, which provides the way to prove the existence of minimizers. This is formalized in the following theorem (see \cite{M90}, \cite{BFK95}, \cite{VG04} for details).

\begin{theo}\label{compCM}
Let $M>0$ be fixed. Then for every $p<+\infty$ the class $C_M$ is compact with respect to the strong topology of the Sobolev space $W^{1,p}_{loc}(\O)$.
\end{theo}

The compactness established in Theorem \ref{compCM} ensures that any minimizing sequence admits a strongly convergent subsequence, which allows to prove the existence of optimal solutions to the general minimization problem \eqref{mpphi}. The precise existence result is as follows.

\begin{theo}\label{exCM}
Let $\phi:\O\times\R\times\R^d\to\overline\R$ be a function such that
\begin{itemize}
\item[(i)]$\phi$ is nonnegative and measurable for the $\sigma$-algebra
$\LL_d\otimes\B\otimes\B_d$, where $\LL_d$ denotes the Lebesgue $\sigma$-algebra and $\B$, $\B_d$ are the Borel $\sigma$-algebras on $\R$ and $\R^d$ respectively;
\item[(ii)]for almost every $x\in\O$ the function $\phi(x,\cdot,\cdot)$ is lower semicontinuous on $\R\times\R^d$.
\end{itemize}
Then, for any fixed $M>0$ the minimization problem \eqref{mpphi} admits at least one solution in $C_M$. In particular, the original Newton's problem \eqref{mpcart} possesses an optimal solution.
\end{theo}

Beyond the constraint on the maximal height, namely the pointwise bound $0\le u\le M$ which ensures the physical and mathematical well-posedness of the problem, one can consider alternative or additional types of constraints imposed on the class of admissible functions. These functions are still required to be nonnegative and concave, reflecting the geometric and dynamical assumptions of the model.

A particularly natural and physically meaningful alternative to the pointwise height constraint is a global bound on the volume of the solid body. In the geometric context considered here, this volume constraint is expressed as
$$\int_\O u\,dx\le V,$$
where $V>0$ is a prescribed upper bound on the admissible volume. Imposing this condition leads to the definition of a new class of admissible functions, denoted by $\V_V$, given by
$$\V_V=\Big\{u\text{ concave, }u\ge0,\ \int_\O u\,dx\le V\Big\}.$$
The concavity assumption allows us to derive a lower bound on the volume in terms of the maximum height of the solid body. Specifically, it is easy to show, by comparing the body itself with the cone of equal height, the inequality
$$\int_\O u\,dx\ge\frac{|\O|}{d+1}\|u\|_{L^\infty}.$$
This estimate immediately implies the inclusion $\V_V\subset C_{(d+1)V/|\O|}$ and allows us to deduce the compactness of the class $\V_V$ via Theorem \ref{compCM}. As a direct consequence, the existence of minimizers in the class $\V_V$ follows from Theorem \ref{exCM}.

Another particularly interesting class of admissible functions arises when one imposes a constraint on the free surface area of the solid body. This constraint translates into the following integral inequality:
$$\int_\O\sqrt{1+|\nabla u|^2}\,dx+\int_{\partial\O}u\,d\HH^{d-1}\le S,$$
where $\HH^{d-1}$ denotes the $d-1$ dimensional Hausdorff measure on the boundary $\partial\O$, and $S>0$ is a prescribed constant that represents the maximal allowed surface area of the body. The first integral accounts for the area of the graph of $u$ over the interior of $\O$, while the second term represents the vertical contribution along the lateral boundary. This leads to the definition of the admissible class $\HH_S$, given by
$$\HH_S=\Big\{u\ge0,\ u\text{ concave, }\int_\O\sqrt{1+|\nabla u|^2}\,dx+\int_{\partial\O}u\,d\HH^{d-1}\le S\Big\}.$$
Once again, the structural assumption that $u$ is concave plays a key role in deriving geometric estimates. In particular, by comparing again the body itself with the cone of equal height, one can establish the inequality
$$S\ge\frac{|\partial\O|}{d}\|u\|_{L^\infty}.$$
As a consequence of this estimate, we obtain the inclusion
$\HH_S\subset C_{Sd/|\partial\O|}$. This allows us to deduce the compactness property for the class $\HH_S$ from those already established for the class $C_M$ (via Theorem \ref{compCM}), and in turn, to infer existence results for minimization problems posed on $\HH_S$ by invoking Theorem \ref{exCM}.

For a more detailed treatment of the Newton's problem in these classes of admissible shapes, we refer the reader to \cite{BW03}, \cite{HKV02}, \cite{W99}.

It is important to emphasize that the Newton minimal resistance problem does not admit smooth (i.e., globally regular) solutions. In fact, one can rigorously establish the following property on the gradient behavior of admissible solutions.

\begin{theo}
Let $u$ be a solution to Newton's minimal resistance problem \eqref{mpcart}. Then, at any point $x\in\O$ where $u$ is differentiable, we have either $\nabla u(x)=0$ or $|\nabla u(x)|\ge1$.
\end{theo}

\begin{rema}
This result immediately precludes the possibility of regular (e.g., $C^1$) solutions. Indeed, consider the upper boundary of the region where the function attains its maximal value, $\{u=M\}$; at points on this boundary, the gradient necessarily satisfies $|\nabla u(x)|\ge1$, then violating any condition of smooth decay. Consequently, any solution must exhibit a certain region of singularity (specifically, discontinuities in $\nabla u$).

At present, it remains an open question whether additional singularities exist elsewhere in the domain and, if so, what their precise structure and location might be. The nature and distribution of such singular points continue to pose a challenging and intriguing problem in the analysis of variational models with geometric constraints as the Newton's minimal resistance problem.
\end{rema}

The original resistance problem was considered by Newton in the framework of {\it radial} functions; more precisely, if $\O$ is a two-dimensional circle of radius $R$ and $u(r)$ only depends on the radial variable, the resistance functional in \eqref{mpcart} takes the form
$$2\pi\int_0^R\frac{r}{1+|u'|^2}\,dr.$$
The corresponding minimization problem
\be\label{radial}
\min\bigg\{\int_0^R\frac{r}{1+|u'|^2}\,dr\ :\ u\text{ concave, }0\le u\le M\bigg\}
\ee
can be solved explicitly by means of the Euler-Lagrange equation
$$ru'=C\big(1+{u'}^2\big)^2\qquad\hbox{on the set }\{u'\ne0\}$$
for a suitable constant $C<0$. In parametric form we obtain $u(r)=M$ on an interval $[0,r_0]$ and
\be\label{ur}
\begin{cases}
\displaystyle r(t)=\frac{r_0}{4t}(1+t^2)^2\\
u(t)=\displaystyle M-\frac{r_0}{4}\Big(-\frac{7}{4}+\frac{3}{4}t^4+t^2
-\ln t\Big)
\end{cases}
\qquad\forall t\in[1,T].
\ee
Here the quantities $r_0$ and $T$ are defined through the strictly increasing function
$$f(t)=\frac{t}{(1+t^2)^2}\left(-\frac{7}{4}+\frac{3}{4}t^4+t^2-\ln
t\right)\qquad\forall t\ge1$$
by setting:
$$T=f^{-1}(M/R),\qquad r_0=\frac{4RT}{(1+T^2)^2}.$$
Notice that $|u'(r)|>1$ for all $r>r_0$ and that $|u'(r^+_0)|=1$; in particular, the derivative $|u'|$ never belongs to the interval $]0,1[$. We emphasize that, for all values of $R$ and $M$, the optimal radial profile invariably exhibits a flat region at its top, a feature originally anticipated by Newton in his foundational analysis.

The explicit expression \eqref{ur} of the optimal radial profile allows to compute the asymptotic behavior of the radius $r_0$ of the upper flat part and of the relative resistance
$$C_0=\frac{2}{R}\int_0^R\frac{r}{1+|u'|^2}\,dr$$
as $M/R\to+\infty$, being $\O$ the disk of radius $R$. We have
$$r_0/R\approx\frac{27}{16}(M/R)^{-3}\;,\qquad C_0\approx\frac{27}{32}(M/R)^{-2}.$$

\begin{figure}[h!]
\centering
{\includegraphics[scale=0.37]{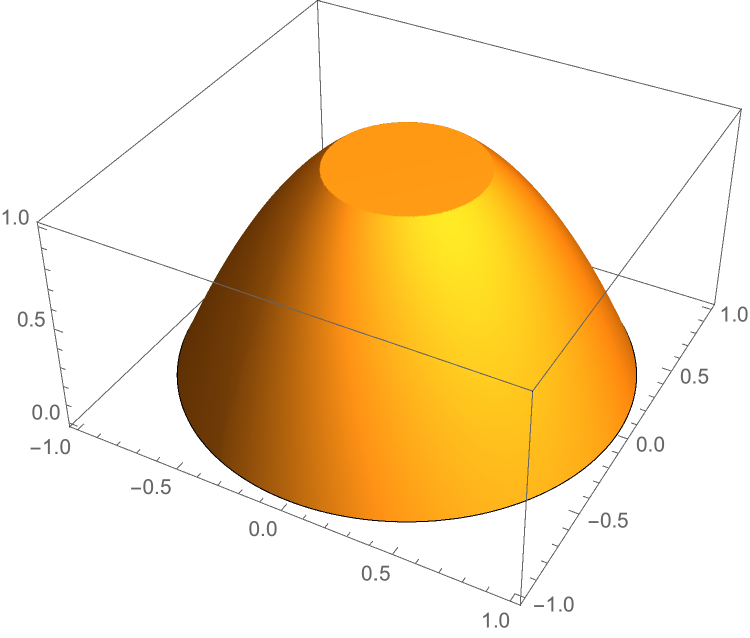}
\includegraphics[scale=0.37]{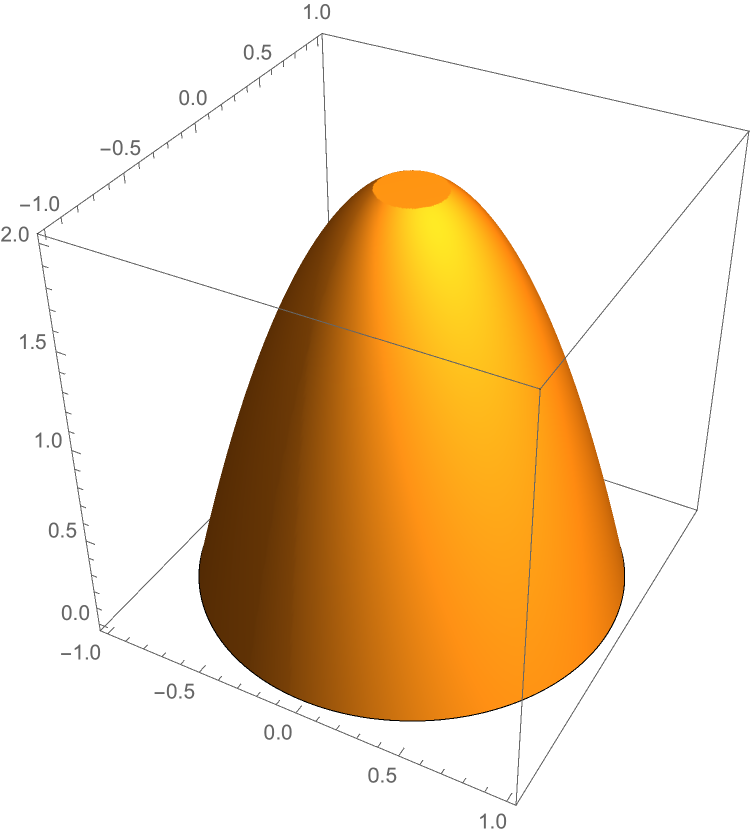}
\includegraphics[scale=0.37]{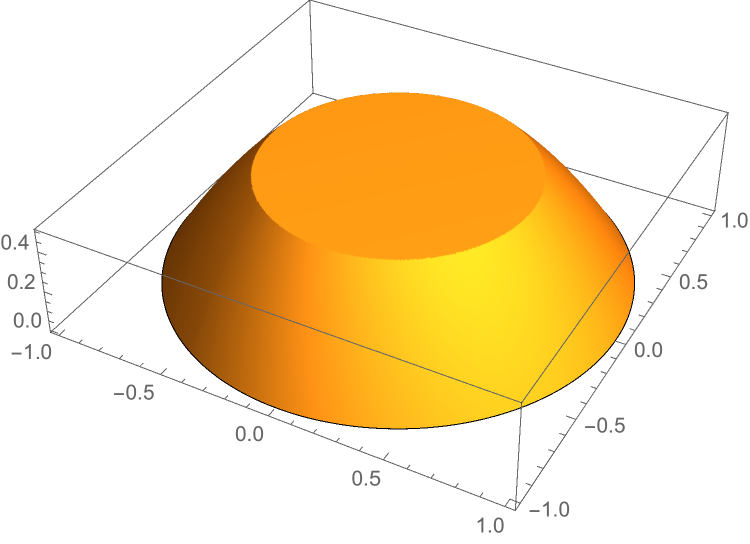}}
\caption{Optimal radial shapes with $R=1$ and $M=1$ (left), $M=2$ (center), $M=0.5$ (right).}\label{fig1}
\end{figure}

A natural and intriguing question that arises in this context is whether, when $\O$ is a circular domain, the optimal solution to Newton's problem of minimal resistance, when restricted to the class of all concave functions, coincides with the radial solution described above. This expectation appears quite reasonable, given the inherent symmetry of the circular domain and the intuitive appeal of radial profiles. However, despite this seemingly natural conjecture, the Newton minimal resistance problem reveals a surprising and nontrivial symmetry-breaking phenomenon. This unexpected behavior is a consequence of the following result, established in \cite{BFK96}.

\begin{theo}\label{flat}
Assume that $u$ is an optimal solution for the Newton problem \eqref{mpcart} which is of class $C^2$ in an open set $\omega\subset\O$ and that $u<M$ in $\omega$. Then we have
\be\label{flatness}
\det\nabla^2 u\equiv0\qquad\hbox{in }\omega.
\ee
\end{theo}

Indeed, while the optimal radial solution given in \eqref{ur} is smooth in the region $\{r>r_0\}$, it fails to satisfy the flatness condition specified in \eqref{flatness} within this domain. This discrepancy has an immediate consequence: it precludes the possibility that the solution, when $\O$ is a circle, remains purely radial.

\begin{figure}[h!]
\centering
\includegraphics[scale=0.7]{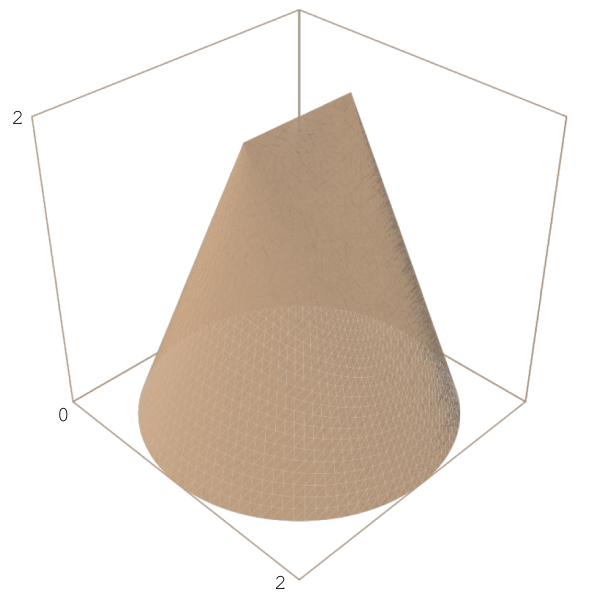}
\caption{Optimal convex shape with $R=1$ and $M=2$.}
\label{fig2}
\end{figure}

Figure \ref{fig2} shows an optimal non-radial shape; we refer the interested reader to \cite{B09}, where numerous images of optimal radial and non-radial bodies are provided.

%%%%%%%%%%%%%%%%%%%%%%%%%%%%%%%%%%%%%%%%%%%%%%%%%%
\section{The intrinsic model}\label{sintri}

An alternative and more geometrically intrinsic formulation of Newton's classical problem of minimal resistance involves expressing the resistance of a solid body $E\subset\R^{d+1}$ as a boundary integral over its surface $\partial E$, without explicit reference to a particular coordinate representation such as its Cartesian graph.

Suppose the convex solid body $E$ is described by the epigraph of a concave function $u:\O\to\R$, where $\O$ is the maximal cross section, so that the boundary $\partial E$ can be represented locally as the graph of $u$. At each point $x\in\partial E$, the unit outward normal vector $\nu(x)$ to the surface is then given by the normalized vector
$$\nu=\frac{(-\nabla u,1)}{\sqrt{1+|\nabla u|^2}}.$$
This representation enables a change of variables from an integral over the domain $\O$ to a surface integral on $\partial E$, yielding the classical expression for the resistance:
\be\label{cube}
\int_\O\frac{1}{1+|\nabla u|^2}\,dx=\int_{\partial E}\big({\nu_n(x)_+}\big)^3\,d\HH^d(x),
\ee
where $\nu_n$ denotes the component of the normal vector in the vertical direction (the $(d+1)$-th coordinate), and $s_+$ denotes the positive part of $s$. This motivates the consideration of a broader class of functionals, given in the general form
\be\label{intrinF}
F(E)=\int_{\partial E}f\big(x,\nu(x)\big)\,d\HH^d(x)
\ee
where $f:\R^{d+1}\times S^d\to\R$ is a given continuous integrand depending on both the position $x$ on the surface and on the unit normal vector $\nu(x)$ at that point. The admissible configurations are taken to be convex bodies $E$ satisfying
$$C_{K,Q}=\Big\{E\text{ convex subset of }\R^{d+1}\ :\ K\subset E\subset Q\Big\},$$
where $K$ and $Q$ are two compact subsets of $\R^{d+1}$ with $K\subset Q$.

In this setting, the Newton minimal resistance problem takes the following intrinsic variational form:
\be\label{intrpb}
\min\bigg\{\int_{\partial E}f\big(x,\nu(x)\big)\,d\HH^d\ :\ E\in C_{K,Q}\bigg\}.
\ee
Again, additional constraints on the volume $|E|$ or on the surface $|\partial E|$ can be imposed. The following result has been obtained in \cite{BG97}.

\begin{theo}\label{maires}
Let $f:\R^{d+1}\times S^d\to\R$ be a bounded continuous function. Then the minimum problem \eqref{intrpb} admits at least one solution.
\end{theo}

\begin{rema}
By taking $K=\overline\O\times\{0\}$ and $Q=\overline\O\times[0,M]$, with $f(x,\nu)=(\nu_n)_+^3$, we recover the classical Newton's problem of minimal resistance.
\end{rema}

%%%%%%%%%%%%%%%%%%%%%%%%%%%%%%%%%%%%%%%%%%%%%%%%%%
\section{The model with positive temperature}\label{stempe}

One of the notable weak points of Newton's classical model lies in its complete disregard for the influence of body regions that lie below the maximal cross-sectional area on the overall resistance of a body $E$. According to this model, only the part above the maximal cross-section contributes to the resistance, while the remaining geometry of the body is treated as irrelevant. This simplification, however, stands in contrast to the experimental evidence, which consistently demonstrate that these sub-maximal regions play a significant role in the dynamics. Accordingly, more refined and realistic models must account for the contribution of the entire geometry of the body, not merely its surface profile above the cross section $\O$.

In Newton's original formulation, the fluid through which the body moves is idealized to consist of particles that remain stationary, with motion attributed solely to the body itself. This abstraction effectively models the fluid as being at zero temperature, where thermal agitation is absent. If, however, one allows the fluid particles to exhibit random, chaotic motion, then it becomes natural to associate such a modification with a physical framework of non-zero temperature environment. In this context, the interaction between the body and the fluid becomes more intricate: the entire boundary $\partial E$ of the body $E$ acquires relevance in the analysis of resistance

The scenario in which the fluid particles are endowed with a positive temperature, allowing for random thermal motion, has been rigorously investigated in the work \cite{PT05}. For the sake of completeness and to facilitate a better understanding of the modifications introduced by this setting, we briefly recall the formulation of that model in what follows.

The following model has been proposed in \cite{PT05} (see also \cite{P12}). The body $E$ is now supposed static, while fluid particles move with a distribution density of velocities $\rho(v)$. After some simple mechanical considerations one finds the expression of the resistance force
$$R(E)=\int_{\partial E}f\big(\nu(x)\big)\,d\HH^d$$
where the vector $f(\nu)$ is given by
$$f(\nu)=-2\nu\int_{\R^d}(v\cdot\nu)_-^2\rho(v)\,dv$$
being $s_-$ the negative part of $s$.

\begin{rema}
For instance, in the classical Newton's case $\rho(v)$ is the Dirac distribution $\delta_{-Vn}$, where $n$ is the vector $(0,0,\dots,1)$, which gives the resistance force in the direction $n$ proportional to the integral
$$\int_{\partial E}\big(\nu_n(x)\big)_+^3\,d\HH^d,$$
which is consistent with the expression \eqref{cube}.
\end{rema}

\begin{figure}[h!]
\centering
\includegraphics[scale=0.7]{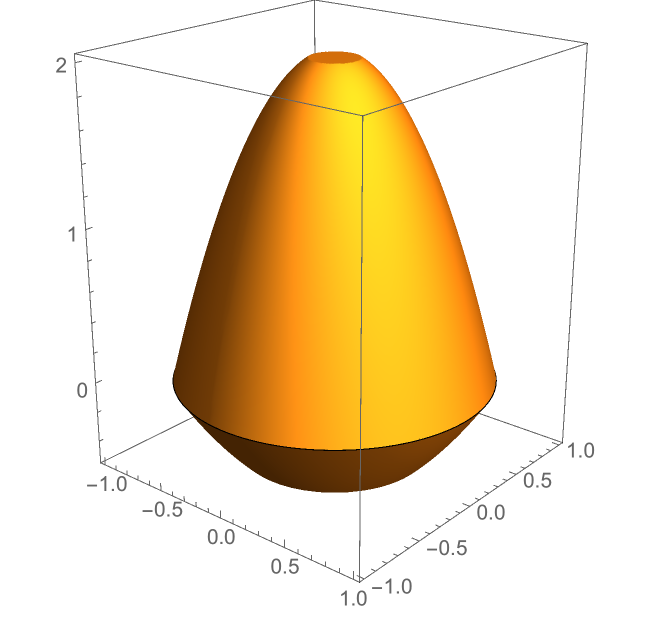}
\caption{Optimal convex radially symmetric shape in the case of positive temperature.}
\label{fig3}
\end{figure}

In Figure \ref{fig3} below a plot of an optimal radially symmetric shape in the case of positive temperature.

%%%%%%%%%%%%%%%%%%%%%%%%%%%%%%%%%%%%%%%%%%%%%%%%%%
\section{Some related problems}\label{srelat}

In recent years, a number of promising and stimulating directions of research have emerged in relation to the classical and modern formulations of problems involving minimal resistance. In what follows, we provide a brief overview of a selection of these novel research trajectories. For a more extended exposition and additional references we refer the interested reader to the monograph \cite{P12}, which offers a detailed presentation of the subject.

\begin{itemize}
\item The Newton's model for resistance is limited to the class of convex bodies, due to the fact that only the first impact between the fluid particles and the surface of the body is taken into account and contributes to the resulting resistance force. This simplification excludes the possibility of multiple reflections or interactions that can occur in more complex geometries. However, the convexity assumption can be relaxed by extending the analysis to a broader class of bodies, namely those for which every fluid particle has at most one impact with the body's surface. This generalized framework preserves the core assumptions of Newton's model allowing a slightly wider variety of shapes. Several investigations in this extended class of bodies can be found in the works \cite{CL01}, \cite{CL02}, \cite{MMOP17}, and \cite{P16}.

\item In situations where fluid particles may have more than a single impact with the surface of the body, the mathematical characterization of resistance becomes significantly more intricate. The simple framework assumed in Newton's model no longer suffices, as one must account for the complex trajectories resulting from multiple reflections. These scenarios introduce a level of dynamical complexity that needs different analytical tools, some of which are related to the theory of {\it billiards}. An interesting attempt to address this more general and challenging setting can be found in the work \cite{P03}, where an analysis of the resistance in the presence of multiple impacts is undertaken.

\item Another limitation in the original formulation of Newton's problem is the absence of tangential friction, that is on the contrary present in real-world fluid-body interactions. In Newton's idealized model, only the normal component of the impact force is considered, neglecting the frictional forces that arise due to the motion between the fluid particles and the surface of the body. This simplification, while mathematically convenient, renders the model less applicable to practical scenarios where such frictional effects play a significant role in determining the overall resistance. To address this physically relevant situation, a modified approach has been proposed in \cite{HKV02}, where the classical framework is extended to incorporate tangential friction into the resistance computation.

\item The case of {\it rotating bodies} introduces yet another interesting and physically significant class of minimal resistance problems. When a body, moving into a medium, has a rotational motion, the expression of the resistance force has to be suitably revised in order to take into account the additional effects induced by rotation. The analysis of resistance in rotating bodies, when in addition tangential friction is present, includes the study of the {\it Magnus effect}, a phenomenon where a spinning body moving through a fluid shows a lateral deflection of its trajectory, well-known in various sports (for instance soccer, tennis, golf, baseball, \dots). A comprehensive treatment of this phenomenon is presented in the monograph \cite{P12}, where an entire chapter is dedicated to its mathematical modeling and physical implications.

\item Another interesting class of problems, closely related to those concerning minimal resistance, involves the concept of {\it low visibility} of bodies. In this context, visibility can be analyzed using mathematical tools analogous to those in the study of resistance, treating the interaction of particles as the deflection of light rays. Of particular interest is the extremal case of {\it invisible bodies}, that leave not deflected almost all incoming rays in a given direction. A rigorous mathematical treatment of these visibility problems, including the construction and analysis of such invisible bodies, can be found in the works \cite{AP09} and \cite{PR11}, to which we refer the interested reader for a detailed presentation.

The opposite behavior occurs when one has problems concerned not with minimizing, but rather maximizing the visibility of an object. In this setting, the objective is to design a body such that every incident light ray is reflected back along a path opposite to its direction of incidence. Such bodies, which possess the remarkable ability to return incoming rays toward their source regardless of the angle of incidence, are commonly referred to as {\it retroreflectors} and play an important role in several applications. For a rigorous mathematical presentation of this class of problems, including the construction of retroreflecting shapes, we refer the interested reader to the article \cite{BKMP11}, \cite{GPT09}, \cite{P11}.

\end{itemize}

%%%%%%%%%%%%%%%%%%%%%%%%%%%%%%%%%%%%%%%%%%%%%%%%%%
\bigskip
\noindent{\bf Acknowledgments. }The work of GB is part of the project 2017TEXA3H {\it``Gradient flows, Optimal Transport and Metric Measure Structures''} funded by the Italian Ministry of Research and University. GB is member of the Gruppo Nazionale per l'Analisi Matematica, la Probabilit\`a e le loro Applicazioni (GNAMPA) of the Istituto Nazionale di Alta Matematica (INdAM).
\bigskip
%%%%%%%%%%%%%%%%%%%%%%%%%%%%%%%%%%%%%%%%%%%%%%%%%%

\bigskip\bigskip\bigskip
\noindent Giuseppe Buttazzo: Dipartimento di Matematica, Universit\`a di Pisa,\\
Largo B. Pontecorvo, 5 -- 56127 Pisa (ITALY)\\
{\tt giuseppe.buttazzo@unipi.it\qquad https://people.dm.unipi.it/buttazzo/}

\end{document}